\documentclass[11pt]{article}

\usepackage{amsmath,amssymb,amsfonts,theorem,makeidx,latexsym,epsfig,subfigure}

\usepackage{color}

\newtheorem{defn}{Definition}[section]

\newtheorem{lemma}[defn]{Lemma}

{\theorembodyfont{\rmfamily}

\newtheorem{ex}[defn]{Example}}

\newtheorem{thm}[defn]{Theorem}

\newtheorem{prop}[defn]{Proposition}

\newtheorem{cor}[defn]{Corollary}

\newtheorem{rem}[defn]{Remark}

\numberwithin{equation}{section}

\newcommand{\ltr}{ L^2(\mathbb R) }

\newcommand{\mn}{\mathbb N}

\newcommand{\mr}{\mathbb R}

\newcommand{\mz}{\mathbb Z}

\newcommand{\mc}{\mathbb C}

\newcommand{\mts}{ \{E_{mb}T_{na}g \}_{m,n \in \mz}}

\def\bp{{\noindent\bf Proof. \ }}

\def\ep{\hfill$\square$\par\bigskip}

\def\bqs{\begin{equation}}

\def\eqs{\tag*{$\square$}\end{equation}\par\bigskip}

\def\la{\langle}

\def\ra{\rangle}

\def\supp{\text{supp}}

\def\vn{\vspace{.1in}\noindent}

\def\bop{\begin{op}\rm}

\def\eop{\end{op}}

\def\bee{\begin{eqnarray}}

\def\ene{\end{eqnarray}}

\def\bes{\begin{eqnarray*}}

\def\ens{\end{eqnarray*}}

\def\bei{\begin{itemize}}

\def\eni{\end{itemize}}

\def\bt{\begin{thm}}

\def\et{\end{thm}}

\def\bc{\begin{cor}}

\def\ec{\end{cor}}

\def\bpr{\begin{prop}}

\def\epr{\end{prop}}

\def\bl{\begin{lemma}}

\def\el{\end{lemma}}

\def\bd{\begin{defn}}

\def\ed{\end{defn}}

\def\bex{\begin{ex}}

\def\enx{\end{ex}}

\def\bfi{\begin{fig}}

\def\efi{\end{fig}}

\newcommand{\nft}{ || f||^2}

\title{On the Gabor frame set for compactly supported continuous functions
}
\date{\today}

\author{Ole Christensen\thanks{Department of Applied Mathematics and Computer Science,
Technical University of Denmark,
Building 303, 2800 Lyngby, Denmark (ochr@dtu.dk)
},
Hong Oh Kim\thanks{
Division of General Studies, UNIST,
50 UNIST-gil, Ulsan 44919,
Republic of Korea (hkim2031@unist.ac.kr)},
Rae Young Kim\thanks{Department of Mathematics,
Yeungnam University,
280 Daehak-Ro, Gyeongsan, Gyeongbuk 38541,
Republic of Korea (rykim@ynu.ac.kr)
}}

\begin{document}

\maketitle

\begin{abstract} We identify a class of continuous compactly supported functions
for which the known part of the Gabor frame set can be extended. At least for functions
with support on an interval of length two, the curve determining the
set touches the known obstructions. Easy verifiable sufficient conditions for
a function to belong to the class are derived, and it is shown that the
B-splines $B_N, N\ge 2,$ as well as certain ``continuous and truncated" versions
of several classical functions (e.g., the Gaussian and the two-sided exponential function) belong to the class. The sufficient conditions for the frame property guarantees the existence of a dual
window with a prescribed size of the support.
\end{abstract}

\begin{minipage}{120mm}
\end{minipage}

\section{Introduction}
Frames is a functional analytic tool to obtain representations of the elements in a Hilbert
space as a (typically infinite) superposition of building blocks. Frames indeed lead to
decompositions that are similar to the ones obtained via orthonormal bases, but with
much greater flexibility, due to the fact that the definition is significantly less
restrictive. For example, in contrast to the case for a basis, the elements in a frame
are not necessarily (linearly) independent, i.e., frames can be redundant.

One of the main manifestations of frame theory is within Gabor analysis, where the
aim is to obtain efficient representations of signals in a way that reflects the time-frequency distribution.
For any $a,b>0,$ consider the translation operator $T_a$ and the modulation operator $E_b,$
both acting on the particular Hilbert space $\ltr,$ given by $T_af(x)=f(x-a),$ respectively $E_bf(x)=e^{2\pi i bx}f(x).$
Given $g\in \ltr,$ the collection of functions $\mts$ is called a {\it (Gabor) frame}
if there exist constants $A,B>0$ such that
\bes A\, \nft \le \sum_{m,n\in \mz} | \la f, E_{mb}T_{na}g\ra|^2 \le B\, \nft, \, \forall f\in \ltr.\ens If at least the upper condition is satisfied, $\mts$ is called a {\it Bessel
sequence.} It is known that for every frame $\mts$ there exists a dual frame
$\{E_{mb}T_{na}h\}_{m,n\in \mz}$ such that each $f\in \ltr$ has the
decomposition
\bee \label{60115a}  f=\sum_{m,n\in \mz} \la f, E_{mb}T_{na}h\ra E_{mb}T_{na}h.\ene
The problem of determining  $g\in \ltr$ and parameters $a,b>0$ such that $\mts$ is
a frame has attracted a lot of attention over the past 25 years.
The {\it frame set} for a function $g\in \ltr$ is defined as the set
\bes \label{4104a} {\cal F}_g:= \left\{(a,b)\in \mr_+^2 \, \big|
\, \mts \, \mbox{is a frame for} \, \ltr\right\}.\ens
Clearly the ``size" of the set ${\cal F}_g$ reflects the flexibility of the
function $g$ in regard of obtaining expansions of the type \eqref{60115a}.
In particular it is known that $ab\le 1$ is necessary for $\mts$ to be a frame and that
the number $(ab)^{-1}$ is a measure of the redundance of the frame; the smaller
the number is, the more redundant the frame will be. Thus a reasonable
function $g$ should lead to a frame $\mts$ for values $(ab)^{-1}$ that are reasonably
close to one. We remark that ${\cal F}_g$ is known to be open if $g$ belongs to the
Feichtinger algebra; see \cite{FK,AFK}.

Until recently the
exact frame set was only known for very few functions:
the Gaussian $g(x)=e^{-x^2}$  \cite{Ly,Se2,SW}, the hyperbolic secant \cite{JS4}, and
the functions $h(x)= e^{-|x|}, k(x)=e^{-x}\chi_{[0, \infty[}(x)$ \cite{Jan9,Jan2003}.
In \cite{GrSt} a characterization was obtained for the class of totally positive functions of
finite type, and based on \cite{Jan03} the frame set for functions $\chi_{[0,c]}, \, c>0,$
was characterized in \cite{DaiSun}.

For the sake of applications of Gabor frames it is essential that
the window $g$ is a continuous function with compact support.  Most of the related
literature deals with special types of functions like truncated
trigonometric functions or various types of splines, see \cite{Del1,Lau,KlSt,KimI}.
Various classes of functions have also been considered, e.g.,
functions yielding a partition of unity \cite{GJ,CKK6}, functions with short
support or a finite number of sign-changes \cite{CKK1,CKK3,CKK8}, or functions that are bounded away from
zero on a specified part of the support \cite{LN}. The case of B-spline generated Gabor systems has attracted special attention, see, e.g., \cite{pr,KlSt,CKK8,LN,GR2015}.

To our best
knowledge the frame set has not been characterized for any
function $g\in C_c(\mr) \setminus \{0\}.$  We will,  among others, consider a class of
functions for which we can extend the known set of parameters
$(a,b)$ yielding a Gabor frame. The class of functions contains
the B-splines $B_N, \, N\ge 2,$ as well as certain ``continuous
and compactly supported variants" of the above functions $g, h$ and other classical functions. Furthermore, the results
guarantees the existence of dual windows with a support size given
in terms of the translation parameter.

In the rest of this introduction we will describe the relevant class of windows and their
frame properties. Proofs of the frame properties are in Section \ref{50910b}, and easy verifiable conditions for a function to belong to the class
are derived in Section \ref{50910a}.

Let us first collect some of the known results concerning frame properties for
continuous compactly supported functions; (i) is classical, and we refer to
\cite{CB} for a proof.

\bpr \label{50910c} Let $N>0$, and assume that $g: \mr \to \mc$ is
a continuous function with $\supp \, g \subseteq [-\frac{N}{2},
\frac{N}{2}]$. Then the following holds: \bei \item[(i)] If $\mts$
is a frame, then $ab< 1$ and $a<N.$ \item[(ii)] {\rm \cite{LN}}
Assume that $0<a<N, 0<b\leq \frac{2}{N+a}$, and
$\inf_{x\in[-\frac{a}{2},\frac{a}{2}]} |g(x)|>0.$ Then   $\{E_{m
b}T_{n a}  g\}_{m,n\in \mz}$ is a frame, and there is a unique
dual window $h\in L^2(\mr)$ such that {\em $\supp \
h=[-\frac{a}{2},\frac{a}{2}]$}. \item[(iii)] {\rm \cite{CKK8}}
Assume that  $\frac{N}{2} \le  a < N $ and $ 0<b< \frac{1}{a}.$ If
$ g(x)>0, \ \ x\in ]-\frac{N}{2}, \frac{N}{2}[,$ then   $\{E_{m
b}T_{n a}  g\}_{m,n\in \mz}$ is a frame. \eni \epr

We will now introduce the window class that will be used in the current paper; it is a subset
of the set of functions $g$ considered in Proposition \ref{50910c} (iii).
The definition is
inspired by certain explicit estimates for B-splines, given by Trebels and Steidl in \cite{TS};
this point will be clear in Proposition \ref{bn-4-1}.
First, fix $N>0$ and $0<a<N$.
Consider the first order difference $\Delta_a f$ and the second order difference $\Delta^2_a f,$ given by
$$\Delta_a f(x)= f(x)-f(x-a), \, \, \, \Delta^2_a f(x)= f(x)-2f(x-a)+f(x-2a).$$
We define the window class as the set of functions
$$V_{N,a}:=\{f \in C(\mr)|\ \supp \
f=[-\frac{N}{2},\frac{N}{2}], \ f\ \text{is real-valued and satisfies (A1)-(A3)} \}, $$
where
\begin{itemize}
    \item[{\rm(A1)}] $f$ is symmetric around the origin;
    \item[{\rm(A2)}] $f$ is strictly increasing on $[-\frac{N}{2},0]$;
    \item[{\rm(A3)}] If $a<\frac{N}{3}$, then $\Delta^2_a f(x) \geq 0, \, x\in [-\frac{N}{2}, -\frac{N}{4}+\frac{3a}{4}]$; if $a\geq\frac{N}{3}$, then
$\Delta^2_a f(x) \geq 0, \, x\in [-\frac{N}{2},0]\cup\{-\frac{N}{4}+\frac{3a}{4}\}$.
\end{itemize}
Note that by the symmetry condition (A1) a function $f\in V_{N,a}$ is completely determined by its
behavior for $x\in [-\frac{N}{2},0]$.
 If $a\geq \frac{N}{3}$,
the point $-\frac{N}{4}+\frac{3a}{4}$
considered in (A3) is not contained in
$[-\frac{N}{2},0]$; however, if desired, the symmetry condition  allows to formulate
the condition $\Delta^2_a f(-\frac{N}{4}+\frac{3a}{4})\geq 0$
alternatively as
\begin{equation} \label{60208f}
   f(\frac{N}{4}-\frac{3a}{4})-2f(-\frac{N}{4}-\frac{a}{4})\geq 0
\end{equation}
because the argument $x-2a$ of the last term in the second order difference is less than $-\frac{N}{2}$.

The definition of $V_{N,a}$ is technical, but we will derive easy verifiable
conditions for a function $g$ to belong to this set in Proposition \ref{bn-4-1}, and also provide
several natural examples of such functions.
Our main result  extends the range of $b>0$
yielding a frame, compared with Proposition \ref{50910c} (ii):

\bt \label{bn-3-1} For $N>0$, let $0<a<N$
and $\frac{2}{N+a}<b\leq \frac{4}{N+3a}$. Assume that $g\in V_{N,a}$. Then the
Gabor system $\{E_{mb}T_{na}g\}_{m,n\in \mz}$ is a frame for $\ltr$, and there is a unique
dual window $h\in L^2(\mr)$ such that   $\supp \ h \subseteq
[-\frac{3a}{2},\frac{3a}{2}]$. \et

\begin{figure}
\begin{center}
\includegraphics[width=4.9in,height=2.1in]{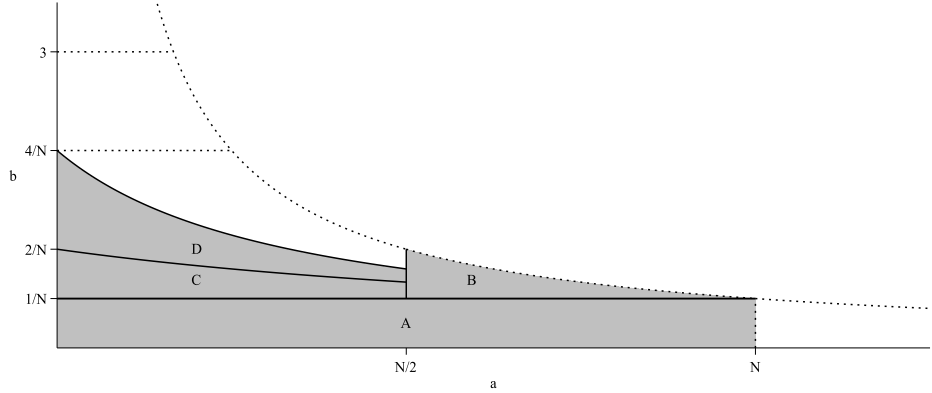}\hfil
\end{center}\caption{The Figure shows the following regions: B as in Proposition
\ref{50910c} (iii), and D as in Theorem \ref{bn-3-1}. The region A
corresponds to the case where the frame operator is a
multiplication operator, and $\mts$ is a frame if $\inf_{x\in
[0,a]} \sum_{n\in \mz} |g(x-na)|^2 >0.$ The region C is a part of the region
determined by Proposition \ref{50910c} (ii), corresponding to
the new findings in the paper \cite{LN}. }
\end{figure}

Membership of a function $g$ in a set $V_{N,a}$ for some $a\in ]0,
N[$ only gives information about the frame properties of $\mts$
for this specific value of the translation parameter $a.$ In order
to get an impression of the frame properties of $\mts$ in a region
in the $(a,b)$-plane, we need to consider a function $g$ that
belongs to $V_{N,a}$ for an interval of $a$-values, preferably for
all $a\in ]0,N[.$  Fortunately several natural functions have this
property. The following list collects some of the results we will
obtain in Section \ref{50910a}.  Considering any $N\in
\mn\setminus\{1 \}$,

\bei
\item The B-spline  $B_N$ of order $N$  belongs to $\bigcap_{0<a<N} V_{N,a}$;

\item The function
$f_N(x):=\cos^{2N-2}(\frac{\pi x}{N})\chi_{[-\frac{N}{2},\frac{N}{2}]}(x)$ belongs to
$\bigcap_{0<a<N} V_{N,a}$;
\item The function
$h_N(x):=\left(e^{-|x|}-e^{-\frac{N}{2}}\right) \chi_{[-\frac{N}{2},\frac{N}{2}]}(x)$
belongs to
$\bigcap_{0<a<N} V_{N,a}$;
\item The function
$g_N(x):=\left(e^{-x^2}-e^{-\frac{N^2}{4}}\right) \chi_{[-\frac{N}{2},\frac{N}{2}]}(x)$
belongs to $\bigcap_{\frac{3N}{7}\leq a<N} V_{N,a}$.
\eni

In particular, Proposition \ref{50910c} and Theorem \ref{bn-3-1} imply that
for $N\in \mn \setminus \{1\}$ the functions $B_N, f_N,$ and $h_N$ generate frames whenever
$0<a<N$ and $0<b\leq \frac{4}{N+3a}$; and
$g_N$ generates a frame whenever
$\frac{3N}{7}\leq a<N$ and $0<b\leq \frac{4}{N+3a}$.

Note that the limit curve $b=\frac{4}{N+3a}$ in Theorem \ref{bn-3-1} touches the
known obstructions for Gabor frames. In fact, for $N=2$ we obtain that
$b\to 2$ whenever $a\to 0$.  Since it is known that the B-spline $B_2$ does not generate a frame
for $b=2$ \cite{Del1,GJ} we can not go beyond this. We also know that
at least for some functions
$g\in \bigcap_{0<a<N} V_{N,a}$ parts of
the region determined by the inequalities
$b<2, a<2, ab <1$
do not belong to the frame set. Considering  for example the B-spline $B_2$,
\cite{LN} shows that the point $(a,b)= (\frac{2}{7},\frac{7}{4})$ does not belong to the frame set.
For $a=\frac{2}{7}$
Theorem \ref{bn-3-1} guarantees the frame property for $b< \frac{7}{5},$ which is close to the obstruction.
These considerations indicate that the frame region in Theorem \ref{bn-3-1} in a quite accurate
way describes the maximally possible frame set below $b=2$ that is valid for all the functions in $V_{N,a}$, at least for $N=2.$

\section{Frame properties for functions $g\in V_{N,a}$} \label{50910b}
The purpose of this section is to prove Theorem \ref{bn-3-1}.
Since the functions $g\in V_{N,a}$ are bounded and have compact support, they generate
Bessel sequences $\mts$ for all $a,b>0.$
By the duality conditions
\cite{RoSh1,Jan2}, two bounded functions $g,h$ with compact support generate dual frames
$\mts$ and $\{E_{mb}T_{na}h\}_{m,n\in \mz}$ for some fixed $a,b>0$ if and only if
\bes \sum_{m\in \mz} g(x-\ell/b+ma)\overline{h(x+ma)} = b\delta_{\ell,0}, \ a.e. \ x\in
[-\frac{a}{2},\frac{a}{2}];\ens in particular, a function $g\in V_{N,a}$
and a bounded real-valued function $h$ with support on $[-\frac{3a}{2},\frac{3a}{2}]$
generate dual Gabor frames
$\{E_{mb}T_{na}g\}_{m,n\in \mz}$ and
$\{E_{mb}T_{na}h\}_{m,n\in \mz}$ for $\ltr$
for some $b\leq \frac{4}{N+3a}$  if and only if the equations
\begin{equation} \label{50911a}
\sum_{m=-1}^{1}
g(x-\ell/b+ma)h(x+ma) = b\delta_{\ell,0}, \ a.e. \ x\in
[-\frac{a}{2},\frac{a}{2}].
\end{equation}
hold for $\ell=0,\pm 1$.
Given $g\in V_{N,a}$ we will therefore  consider the $3\times 3$ matrix-valued function
$G$ on $[-\frac{a}{2},\frac{a}{2}]$ defined by
\bes G(x):=\begin{pmatrix} g(x-\frac{\ell}{b}+ma)
\end{pmatrix}_{-1\leq \ell,m\leq 1}=
\begin{pmatrix}
g(x+\frac{1}{b}-a) & g(x+\frac{1}{b}) & g(x+\frac{1}{b}+a)\\
g(x-a) & g(x) & g(x+a)\\
g(x-\frac{1}{b}-a) & g(x-\frac{1}{b}) & g(x-\frac{1}{b}+a)
\end{pmatrix}.
\ens
In terms of the $G(x)$ the condition \eqref{50911a} simply means that
\begin{equation} \label{bn-15}
  G(x)
 \begin{pmatrix}
h(x-a) \\
h(x)\\
h(x+a)
\end{pmatrix}
=
\begin{pmatrix}
0 \\
b\\
0
\end{pmatrix}, \ a.e.  \ x\in [-\frac{a}{2},\frac{a}{2}].
\end{equation}

We will show that the matrix $G(x)$ is invertible for all $x\in [-\frac{a}{2},\frac{a}{2}]$; this will
ultimately give us a bounded and compactly supported function $h$ satisfying \eqref{50911a}
and hereby prove Theorem \ref{bn-3-1}.
The invertibility of $G(x)$ will be derived as a consequence of a series
of lemmas, where we first consider $x\in [-\frac{a}{2},0].$ Note that the proof of the first result does not use
the property (A3):

\bl \label{bn-11-1}
For $N>0$,
let $0<a<N$ and $\frac{2}{N+a}< b\leq \frac{4}{N+3a}$.
Assume that $g\in V_{N,a},$ and let $x\in[-\frac{a}{2},0].$
Then the following hold:
\begin{itemize}
    \item[{\rm (a)}] $g(x+\frac{1}{b}+a) \leq g(x+\frac{1}{b}) < g(x+\frac{1}{b}-a)\neq 0;$

    \item[{\rm (b)}] $g(x-\frac{1}{b}-a) \leq g(x-\frac{1}{b}) < g(x-\frac{1}{b}+a)\neq 0;$

    \item[{\rm (c)}] $g(x)>g(x-a)$ and
                     $g(x)\geq g(x+a),$ with equality only for $x=-\frac{a}{2}$.
\end{itemize}
\el

\bp For (a), let $x\in [-\frac{a}{2},0]$. Using
 $b\le
\frac{4}{N+3a}$ and $a<N$,
\bes  x+\frac{1}{b}-a \ge  \frac{1}{b}-\frac{3a}{2} \geq
\frac{N+3a}{4}-\frac{3a}{2}  >-\frac{N}{2}.\ens
It follows that \bee && x+\frac{1}{b}-a \subseteq
[\frac{1}{b}-\frac{3a}{2},\frac{1}{b}-a] \subseteq
]-\frac{N}{2},\frac{N}{2}[ \label{bn-18}.\ene Using that $b\le
\frac{4}{N+3a} < \frac{1}{a}$ it follows that $-2x- \frac{2}{b}+a \le 2
(a-  \frac{1}{b})<0;$ thus \bee -x-\frac{1}{b}+a < x+\frac{1}{b}.
\label{bn-21} \ene Since $g(x) >0$ for $x\in
]-\frac{N}{2},\frac{N}{2}[$, we have by \eqref{bn-18} that
$g(x+\frac{1}{b}-a)\neq 0.$ By (A2) we know that $g$ is strictly decreasing
on $[0,\frac{N}{2}]$. If $x+\frac{1}{b}-a\leq 0$, then we have by
\eqref{bn-18}, \eqref{bn-21} and the symmetry of $g$ that
$$g(x+\frac{1}{b}-a)=g(-x-\frac{1}{b}+a)>g(x+\frac{1}{b})\geq g(x+\frac{1}{b}+a);$$
if $x+\frac{1}{b}-a > 0$, then we have
$ g(x+\frac{1}{b}-a)>g(x+\frac{1}{b})\geq g(x+\frac{1}{b}+a).$
Hence (a) holds. Similarly, (b) and (c) hold.
\ep

We now show that if $g\in V_{N,a}$ and $a\ge N/3,$ the condition (A3) automatically holds on
a larger interval.

\begin{lemma} \label{bn-45}
For $N>0$,  let $0<a<N$. Assume that $g\in V_{N,a}$. Then
$\Delta^2_a g(x) \geq 0, \, \forall x\in [-\frac{N}{2}, -\frac{N}{4}+\frac{3a}{4}]$.
\end{lemma}
\bp  It suffices to show that for $a\geq \frac{N}{3}$,
$$\Delta^2_a g(x) \geq 0, \, \forall x\in [0, -\frac{N}{4}+\frac{3a}{4}].$$
We first note that (A1) and (A2) imply that for $x\in[0,-\frac{N}{4}+\frac{3a}{4}]$,
\begin{equation}\label{bn-43}
g(x)\geq g(-\frac{N}{4}+\frac{3a}{4})\text{ and }  g(x-a)\leq g(-\frac{N}{4}-\frac{a}{4}).
\end{equation}
For $a\geq \frac{N}{3}$,  we have $-\frac{N}{4}-\frac{5a}{4}<-\frac{N}{2}$;
due to the compact support of $g$, $f(-\frac{N}{4}-\frac{5a}{4})=0$;  thus
using (A1) again,
\begin{eqnarray*}
\Delta^2_a g(-\frac{N}{4}+\frac{3a}{4})&=&
g(-\frac{N}{4}+\frac{3a}{4})-2g(-\frac{N}{4}-\frac{a}{4})\\
&=&
g(\frac{N}{4}-\frac{3a}{4})-2g(-\frac{N}{4}-\frac{a}{4}).
\end{eqnarray*}
Together with \eqref{bn-43} this shows that
$$0 \leq  \Delta^2_a g(-\frac{N}{4}+\frac{3a}{4}) \leq
\Delta^2_a g(x), \ \forall x\in[0,-\frac{N}{4}+\frac{3a}{4}],  $$
as desired.  \ep

Let $G_{ij}(x)$ denote the $ij$-th minor of $G(x)$, the determinant of the submatrix obtained
by removing the $i$-th row and the $j$-th
column from $G(x)$.

\bl \label{bn-12-1} For $N>0$,  let $0<a<N$ and
$\frac{2}{N+a}< b\leq \frac{4}{N+3a}$. Assume that $g\in V_{N,a}$ and let $x\in[-\frac{a}{2},0].$ Then the following hold:
\begin{itemize}
    \item[{\rm (a)}] $G_{21}(x)\geq 0$, and equality holds
    iff $g(x+\frac{1}{b})=g(x+\frac{1}{b}+a)=0;$
    \item[{\rm (b)}] $G_{23}(x)\geq 0$, and equality holds
    iff $g(x-\frac{1}{b}-a)=g(x-\frac{1}{b})=0;$
    \item[{\rm (c)}]  $G_{22}(x) \geq  G_{21}(x) +G_{23}(x)$.
\end{itemize}
 \el
\bp Since $g\ge 0,$ (a) and (b) follow from Lemma \ref{bn-11-1} (a) \& (b).
For
(c), we note that
$\Delta_a g(x)= g(x)-g(x-a)=g(-x)-g(-x+a)=-\Delta_a g(-x+a)$ by the symmetry of $g$. Now a direct
calculation shows that
\begin{eqnarray*}
&& G_{22}(x) - G_{21}(x) - G_{23}(x) \\
&&=-\Delta_a g( x+\frac{1}{b}) \Delta_a g (x-\frac{1}{b}+a)
+\Delta_a g(x+\frac{1}{b}+a) \Delta_a g(x-\frac{1}{b}) \\
&&=
\Delta_a g(-x-\frac{1}{b}+a) \Delta_a g(x-\frac{1}{b}+a)
- \Delta_a g(-x-\frac{1}{b}) \Delta_a g(x-\frac{1}{b}).
\end{eqnarray*}
Hence it suffices to show that
\begin{equation*} 
\Delta_a g (x-\frac{1}{b}+a) \geq \Delta_a g(x-\frac{1}{b}) , \ x\in [-\frac{a}{2},0]
\end{equation*}
and
\begin{equation*} 
\Delta_a g(-x-\frac{1}{b}+a)\geq \Delta_a g(-x-\frac{1}{b}), \ x\in [-\frac{a}{2},0],
\end{equation*}
or, equivalently, that
$\Delta^2_a g (x-\frac{1}{b}+a) \geq 0$
and
$\Delta^2_a g(-x-\frac{1}{b}+a)\geq 0$, both for $x\in [-\frac{a}{2},0].$ Since
$[-\frac{1}{b}+\frac{a}{2},-\frac{1}{b}+a]\cup [-\frac{1}{b}+a,-\frac{1}{b}+\frac{3a}{2}]=
[-\frac{1}{b}+\frac{a}{2},-\frac{1}{b}+\frac{3a}{2}],$
this means precisely that
\begin{equation} \label{bn-31}
\Delta^2_a g (x) \geq 0, \ x\in
[-\frac{1}{b}+\frac{a}{2},-\frac{1}{b}+\frac{3a}{2}].
\end{equation}
We note that for $\frac{2}{N+a}< b\leq \frac{4}{N+3a}$,
$$-\frac{N}{2}\leq -\frac{1}{b}+\frac{a}{2}, \ -\frac{1}{b}+\frac{3a}{2}\leq -\frac{N}{4}+\frac{3a}{4}.$$
Together with Lemma \ref{bn-45} this implies that \eqref{bn-31}   holds, as
desired.
\ep

After this preparation we can now show that $G(x)$ is indeed invertible for
$x\in [- \frac{a}{2},\frac{a}{2}]$ under the assumptions in Theorem \ref{bn-3-1}.

\bc \label{bn-14-1}

For $N>0$,  let $0<a<N$ and $\frac{2}{N+a}<
b\leq \frac{4}{N+3a}$. Assume that $g\in V_{N,a}$. Then $\det
G(x)\neq 0$ for $x\in [-\frac{a}{2},\frac{a}{2}]$.

\ec

\bp First, consider $x\in [-\frac{a}{2}, 0]$.
By Lemma \ref{bn-12-1} (c), we have
\begin{eqnarray*}
\det G(x)&=&-g(x-a) G_{21}(x)+g(x)G_{22}(x)-g(x+a)G_{23}(x) \\
&\geq& -g(x-a) G_{21}(x)+g(x)(G_{21}(x)+G_{23}(x))-g(x+a)G_{23}(x) \\
&=& (g(x)-g(x-a)) G_{21}(x)+(g(x)-g(x+a))G_{23}(x)\\
 &=:& A_N(x).
\end{eqnarray*}
Using Lemma \ref{bn-11-1} (c) and Lemma \ref{bn-12-1} (a) \& (b),
$$(g(x)-g(x-a)) G_{21}(x) \geq 0 \ \ \mbox{and} \ \
(g(x)-g(x+a))G_{23}(x) \geq 0.$$
Thus
$A_N(x)\geq 0, \
x\in [-\frac{a}{2},0].$ If $A_N(x) > 0$ for all
$x\in [-\frac{a}{2},0]$  the proof is completed; thus the rest of the proof will focus on
the case where
$A_N(x_0)=0$ for some
$x_0\in[-\frac{a}{2},0]$. In this case Lemma \ref{bn-11-1} (c) shows that
either
\bee \label{51117b} G_{21}(x_0)=0,\ G_{23}(x_0)=0\ene
or \bee \label{51117a} G_{21}(x_0)=0,\ x_0=-\frac{a}{2}.\ene
The case \eqref{51117a} actually can not occur. Indeed, if  $G_{21}(-\frac{a}{2})=0$, then by
Lemma \ref{bn-12-1} (a), we have $g(-\frac{a}{2}+\frac{1}{b})=0;$
by the symmetry of $g$ this would imply that
$g(\frac{a}{2}-\frac{1}{b})=0,$ which contradicts
Lemma \ref{bn-11-1} (b) with $x=-\frac{a}{2}$.
Thus we only have to deal with the case \eqref{51117b}.
If $G_{21}(x_0)= G_{23}(x_0)=0$,
then by Lemma \ref{bn-12-1} (a) \& (b), we have
$$g(x_0+\frac{1}{b})=g(x_0+\frac{1}{b}+a)
=g(x_0-\frac{1}{b}-a)=g(x_0-\frac{1}{b})=0.$$ Inserting this information into
the entries of the matrix $G(x_0)$ and applying
Lemma \ref{bn-11-1} yields that
$$\det G(x_0)= g(x_0+\frac{1}{b}-a)g(x_0)
g(x_0-\frac{1}{b}+a)> 0,$$ as desired.
This completes the proof that $G(x) >0$ for $x\in [-\frac{a}{2},0]$.
Since $g$ is symmetric around the origin, we have
\begin{eqnarray*}
  \det G(-x)
  &=&
  \det  \begin{pmatrix}
g(-x-\frac{\ell}{b}+ma)
\end{pmatrix}_{-1\leq \ell,m\leq 1}
=
  \det  \begin{pmatrix}
g(x+\frac{\ell}{b}-ma)
\end{pmatrix}_{-1\leq \ell,m\leq 1} \\
&=&
  - \det  \begin{pmatrix}
g(x-\frac{\ell}{b}-ma)
\end{pmatrix}_{-1\leq \ell,m\leq 1}
=  \det  \begin{pmatrix}
g(x-\frac{\ell}{b}+ma)
\end{pmatrix}_{-1\leq \ell,m\leq 1} \\
&=& \det G(x).
\end{eqnarray*} Thus $G(x)$ is also invertible for $x\in ]0, \frac{a}{2}]$.
\ep

We are now ready to prove Theorem \ref{bn-3-1}.

\vn{\bf Proof of Theorem \ref{bn-3-1}:}  By Corollary \ref{bn-14-1} and
continuity of $g$, $\inf_{x\in [-\frac{a}{2},\frac{a}{2}]} |\det
G(x)|>0$. We define $h$ on $[-\frac{3a}{2},\frac{3a}{2}]$ by
$$
 \begin{pmatrix}
h(x-a) \\
h(x)\\
h(x+a)
\end{pmatrix}
= G ^{-1}(x)
\begin{pmatrix}
0 \\
b\\
0
\end{pmatrix}, \ x\in [-\frac{a}{2},\frac{a}{2}],$$
which is a bounded function. On $\mr\setminus
[-\frac{3a}{2},\frac{3a}{2}]$, put $h(x)=0$. It follows immediately by definition of $h$ that
then $g$ and $h$ are dual windows.
\ep

\begin{rem}
{\em We note that the above approach is tailored to the
region of parameters $(a,b)$ in Theorem \ref{bn-3-1}. For example, it does not apply to the
region considered in Proposition \ref{50910c} (ii). In fact, if $0<b\leq \frac{2}{N+a}$, then
the first row of $G(x)$ for $x\in
[\frac{N}{2}-\frac{1}{b}+a,\frac{a}{2}]$ is the zero vector; and the
third row of $G(x)$ for $x\in
[-\frac{a}{2},\frac{1}{b}-a-\frac{N}{2}]$ is the zero vector. Hence
we have $\inf_{x\in [-\frac{a}{2},\frac{a}{2}]} \det G(x)=0$. }
\end{rem}

The conditions for $g\in V_{N,a}$ are technical.  We will now give an example, showing that
the conclusion about the size of the support of the dual window in Theorem \ref{bn-3-1}
might break down if $g\notin V_{N,a}.$

\begin{ex} Let us consider the parameters $N=5, \, a=1.$ Consider a symmetric and continuous function $g$ on $\mr$ with $\supp
\ g=[-\frac{5}{2},\frac{5}{2}]$; assume further that $g$ is increasing on
$[-\frac{N}{2},0]=[-\frac{5}{2},0],$
and that
$$g(0)=12, \ g(-1)=10, \ g(-\frac{4}{3})= 5, \ g(-\frac{7}{3})=3.$$
Then $\Delta_a^2 g(-\frac{4}{3})=g(-\frac{4}{3})-2g(-\frac{7}{3})+g(-\frac{10}{3})=-1<0.$
Hence $g$ does not satisfy (A3) at $x=-\frac{4}{3}$, i.e. $g\notin V_{N,a}.$

We will show that for  $b=\frac{3}{N+2a}=\frac{3}{7}$ there does not exist a bounded real-valued function $h\in L^2(\mr)$ with $\supp \ h \subseteq [-\frac{3a}{2},\frac{3a}{2}]$, such that the duality conditions
\eqref{50911a} hold. In order to obtain a contradiction, let us assume that such a dual
window indeed exists.
Let
$x_0=0$. Then
$$ x_0-a=-x_0-a = -1, \ x_0-\frac{1}{b}+a=-x_0-\frac{1}{b}+a=-\frac{4}{3}, $$
$$ x_0-\frac{1}{b}=-x_0-\frac{1}{b}=-\frac{7}{3},\ x_0-\frac{1}{b}-a= -x_0-\frac{1}{b}-a =-\frac{10}{3},$$  and consequently
\bes G(x_0)=
\begin{pmatrix}
g(x_0+\frac{1}{b}-a) & g(x_0+\frac{1}{b}) & g(x_0+\frac{1}{b}+a)\\
g(x_0-a) & g(x_0) & g(x_0+a)\\
g(x_0-\frac{1}{b}-a) & g(x_0-\frac{1}{b}) & g(x_0-\frac{1}{b}+a)
\end{pmatrix}
=
\begin{pmatrix}
5 & 3 & 0\\
10 & 12 & 10\\
0 & 3 & 5
\end{pmatrix}.
\ens
By the continuity  of $g$, there exist
continuous functions $\epsilon_{ij}(x)$, $1\leq i,j\leq 3$ such that
$$G(x)=
\begin{pmatrix}
5+  \epsilon_{11}(x)  & 3+ \epsilon_{12}(x) & \epsilon_{13}(x)\\
10+\epsilon_{21}(x) & 12+\epsilon_{22}(x) & 10+\epsilon_{23}(x)\\
\epsilon_{31}(x) & 3+\epsilon_{32}(x) & 5+\epsilon_{33}(x)
\end{pmatrix}
$$
and
$$\epsilon_{ij}(x) \rightarrow 0 \text{ as } x\rightarrow x_0.$$
Then \eqref{bn-15} implies that
\begin{equation*}
\begin{pmatrix}
5+  \epsilon_{11}(x)  & 3+ \epsilon_{12}(x) & \epsilon_{13}(x) \\
10+\epsilon_{21}(x) & 12+\epsilon_{22}(x) & 10+\epsilon_{23}(x)\\
\epsilon_{31}(x) & 3+\epsilon_{32}(x) & 5+\epsilon_{33}(x)
\end{pmatrix}
 \begin{pmatrix}
h(x-a) \\
h(x)\\
h(x+a)
\end{pmatrix}
=
\begin{pmatrix}
0 \\
b\\
0
\end{pmatrix}
\end{equation*} for a.e. $x\in [-\frac{a}{2}, \frac{a}{2}].$
By elementary row operations, this leads to
\begin{equation*}
\begin{pmatrix}
5+  \epsilon_{11}(x)  & 3+ \epsilon_{12}(x) & \epsilon_{13}(x) \\
\eta_1(x) & \eta_2(x) & \eta_3(x)\\
\epsilon_{31}(x) & 3+\epsilon_{32}(x) & 5+\epsilon_{33}(x)
\end{pmatrix}
 \begin{pmatrix}
h(x-a) \\
h(x)\\
h(x+a)
\end{pmatrix}
=
\begin{pmatrix}
0 \\
b\\
0
\end{pmatrix}, \, a.e. \, x\in [-\frac{a}{2}, \frac{a}{2}],
\end{equation*}
where $\eta_i(x):=\epsilon_{2i}(x)-2\epsilon_{1i}(x)-2\epsilon_{3i}(x)$ for $1\leq i \leq 3$.
Ignoring a possible set of measure zero and using that $h$ is a bounded function,  this  implies that
\begin{eqnarray*}
b=\eta_1(x)h(x-a)+\eta_2(x)h(x) + \eta_3(x)h(x+a) \rightarrow  0
\end{eqnarray*}
as $x\rightarrow x_0$.
This is a contradiction.
\ep
\end{ex}

On the other hand, the condition $g\in V_{N,a}$ is not a necessary
condition for $\{E_{mb}T_{na}g \}_{m,n\in\mz}$ to be a frame in
the considered region. For example, let $N=1$, $a=\frac{1}{4},
b=2$ and take $g_1(x):=\left(e^{-x^2}-e^{-\frac{1}{4}}\right)
\chi_{[-\frac{1}{2},\frac{1}{2}]}(x).$
Then elementary calculations show that
(A3) does not
hold for $x\in [-\frac{1}{10},-\frac{1}{16}] \subset
[-\frac{N}{2},-\frac{N}{4}+\frac{3a}{4}]$. But since
$\det
G(x)>0$ for $x\in [-\frac{a}{2},\frac{a}{2}]$,
 one can prove
 that $\{E_{mb}T_{na}g_1
\}_{m,n\in\mz}$ is a frame by following the steps
in the proof of Theorem \ref{bn-3-1}.

\section{The set $V_{N,a}$} \label{50910a}

In this section we give easy verifiable sufficient conditions for a function $g$
to belong to $V_{N,a}$.
Recall that a continuous function $f: \mr \to \mc$
with $\supp \, f= [-\frac{N}{2}, \frac{N}{2}]$ is piecewise continuously differentiable  if there exist  finitely many $x_0=-\frac{N}{2} < x_1 < \cdots < x_n=\frac{N}{2}$ such that
    \begin{itemize}
        \item[(1)] $f$ is continuously differentiable on $]-x_{i-1},x_i[$ for every $i \in \{1,\cdots, n\}$;
        \item[(2)] the one-sided limits $\lim_{x\rightarrow x_{i-1}^+} f^\prime (x)$ and
        $\lim_{x\rightarrow x_{i}^-} f^\prime (x)$ exist for every $i \in \{1,\cdots, n\}$.
    \end{itemize}
Note that if $g$ is a continuous and piecewise continuously differentiable function,
the fundamental theorem of calculus yields that
\begin{equation}\label{bn-10}
\Delta_a g(x)=\int_{x-a}^x g^\prime(t)dt \text{  for all } x\in  \mr, \, a>0.
\end{equation}
In order to avoid a tedious presentation, we will forego to mention the points
where a piecewise continuously differentiable function is not differentiable, e.g., in
conditions (c) and (d) in the
following Proposition \ref{bn-4-1}. The result is inspired by explicit calculations for B-splines,
due to Trebels and Steidl; see Lemma 1 in \cite{TS}.

\bpr \label{bn-4-1}
Let $N>0$ and
assume that a continuous and piecewise continuously differentiable function $g: \mr \to \mr$  with $\supp\ g =[-\frac{N}{2},\frac{N}{2}]$ satisfies the following conditions:
\begin{itemize}
    \item[{\rm (a)}] $g$ is symmetric around the origin;
    \item[{\rm (b)}] $g$ is strictly increasing on $[-\frac{N}{2},0]$;
    \item[{\rm (c)}] $g^\prime$ is  increasing on $]-\frac{N}{2},-\frac{N}{4}]$;
    \item[{\rm (d)}]  $g^\prime (-x-\frac{N}{2}) \leq g^\prime (x)$ for $x\in [-\frac{N}{4},0[$.
\end{itemize}
Then $g\in \cap_{0<a<N}V_{N,a}$. \epr
\bp Note that the conditions
(a) and (b) are exactly the same as (A1) and (A2). Thus we will prove (A3). In the
entire argument we will assume that $g$ is differentiable; an elementary consideration
then extends the result to the case of piecewise differentiable functions.
Let us first consider $x\le - \frac{N}{4}.$ Then by the mean value theorem
\bes \Delta_a^2g(x) & = & a \left( \frac{g(x)-g(x-a)}{a} - \frac{g(x-a)-g(x-2a)}{a}\right)
=  a (g^\prime(\xi)-g^\prime(\eta))\ens
for some $\eta\in [x-2a, x-a], \, \xi \in [x-a,x].$ Since $g^\prime$ is increasing
up to $x=-\frac{N}{4}$ this proves that $\Delta_a^2g(x)\ge 0$ whenever
$x\le - \frac{N}{4}.$

We now consider $x\in [-\frac{N}{4},\min(-\frac{N}{4}+\frac{3a}{4},0)].$ Then
$$x-a \le \min(-\frac{N}{4}- \frac{a}{4},-a)\le -\frac{N}{4}.$$
Then
\bee \notag \Delta_a^2g(x)
& = & \notag
\int_{x-a}^{-\frac{N}{4}} (g^\prime(t)-g^\prime(t-a))\,dt+\int_{-\frac{N}{4}}^x (g^\prime(t)-g^\prime(t-a))\,dt
\\ \notag  & \ge &  \int_{x-a}^{-\frac{N}{4}} (g^\prime(t)-g^\prime(t-a))\,dt+\int_{-\frac{N}{4}}^x (g^\prime(-t-\frac{N}{2})-g^\prime(t-a))\,dt  \\ \notag  & = &
g(-\frac{N}{4})-g(-\frac{N}{4}-a)-g(x-a)+g(x-2a) \\ \notag & \ & - g(-x- \frac{N}{2})-g(x-a)+g(\frac{N}{4}-\frac{N}{2})+g(-\frac{N}{4}-a)
\\ \label{60207a} & = & g(-\frac{N}{4})-g(x-a) - [g(-a-\frac{N}{4})-g(x-2a)]  \\ \label{60207b} & \ & + g(-\frac{N}{4})-g(x-a) -[g(-x-\frac{N}{2})-g(-a-\frac{N}{4})].
\ene
We now consider the terms \eqref{60207a} and \eqref{60207b} separately.
For
\eqref{60207a}, by the mean value theorem,
\bes & \ & g(-\frac{N}{4})-g(x-a) - [g(-a-\frac{N}{4})-g(x-2a)] \\ & = & (-\frac{N}{4}-x+a) \left(\frac{g(-\frac{N}{4})-g(x-a)}{-\frac{N}{4}-x+a} - \frac{g(-a-\frac{N}{4})-g(x-2a)}{-\frac{N}{4}-x+a}\right)
\\ & = &  (-\frac{N}{4}-x+a)\left(g^\prime(\xi)- g^\prime(\eta)\right)\ens
for some $\eta \in [x-2a, -a-\frac{N}{4}], \, \xi \in [x-a,-\frac{N}{4}].$ Since $-a-\frac{N}{4} \le x-a$
we have $\eta \le \xi \le -\frac{N}{4};$ thus, by assumption (c), $g^\prime(\eta) \le g^\prime(\xi).$
Recalling that $x\le - \frac{N}{4}+\frac{3a}{4}< - \frac{N}{4}+a$ we have$-\frac{N}{4}-x+a >0;$ thus, we conclude
that the term in \eqref{60207a} indeed is nonnegative.

For \eqref{60207b} we split into two cases.  If $x-a \ge -x-\frac{N}{2},$ exactly the same argument
as for \eqref{60207a} works. If $x-a< -x-\frac{N}{2}$ we perform the same argument after a
rearrangement of the terms. Indeed,
\bes & \ & g(-\frac{N}{4})-g(x-a) -[g(-x-\frac{N}{2})-g(-a-\frac{N}{4})]   \\ & = & g(-\frac{N}{4})- g(-x-\frac{N}{2})- [g(x-a)-g(-a-\frac{N}{4})] \\ & = &  (x+\frac{N}{4}) \left(\frac{g(-\frac{N}{4})- g(-x-\frac{N}{2})}{x+\frac{N}{4}}- \frac{g(x-a)-g(-a-\frac{N}{4})}{x+\frac{N}{4}}\right) \\ & = & (x+\frac{N}{4}) \left(g^\prime(\xi)- g^\prime(\eta)\right)\ens
for some $\eta \in [-a-\frac{N}{4}, x-a], \, \xi \in [-x-\frac{N}{2},-\frac{N}{4}].$ As before this implies
that \eqref{60207b} is nonnegative.

We now consider $x=- \frac{N}{4}+\frac{3a}{4};$
according to \eqref{60208f} we must prove that
\bee \label{60208b} g(\frac{N}{4}-\frac{3a}{4}) - 2g(- \frac{N}{4}-\frac{a}{4})\ge 0.\ene

First, if $a > \frac{2N}{3}$ (i.e., $- \frac{N}{4}+\frac{3a}{4} > \frac{N}{4}$) we have
\bes & \ & g(\frac{N}{4}-\frac{3a}{4}) - 2g(- \frac{N}{4}-\frac{a}{4})
\\ & = & \frac{N-a}{2} \left( \frac{g(\frac{N}{4}-\frac{3a}{4}) - g(-\frac{N}{4}-\frac{a}{4} )}{\frac{N-a}{2}}-
\frac{g(- \frac{N}{4}-\frac{a}{4})-g(\frac{-3N}{4}+\frac{a}{4})}{\frac{N-a}{2}} \right)
\\ & = &  \frac{N-a}{2}( g^\prime(\xi)-g^\prime(\eta))\ens
for some $\xi \in [-\frac{N}{4}-\frac{a}{4},\frac{N}{4}-\frac{3a}{4}],
\eta \in [-\frac{3N}{4}+\frac{a}{4}, -\frac{N}{4}-\frac{a}{4}];$ thus \eqref{60208b} holds by assumption (c).

Now assume that $\frac{N}{2} < a \le \frac{2N}{3}$; then
$- \frac{N}{4}+\frac{3a}{4} \le \frac{N}{4},$ and
\bes - \frac{N}{4} \le \frac{N}{4}- \frac{3a}{4} \le 0, \, \, - \frac{N}{4}-\frac{a}{4}
< - \frac{N}{2} + \frac{a}{4} <0.\ens
Thus,
\bes & \ & g(\frac{N}{4}-\frac{3a}{4}) - 2g(- \frac{N}{4}-\frac{a}{4})\\
& \ge & g(-\frac{N}{4})-g(- \frac{N}{4}-\frac{a}{4})- (g(-\frac{N}{4}-\frac{a}{4} )-g(-\frac{N}{4}-\frac{2a}{4} ));\ens
this can again be expressed in terms of a difference $g^\prime(\xi)-g^\prime(\eta)$
with $\xi \in [- \frac{N}{4}-\frac{a}{4}, - \frac{N}{4}],
\eta \in [-\frac{N}{4}-\frac{2a}{4}, -\frac{N}{4}-\frac{a}{4}]$ and is hence positive.

Finally, we assume that $\frac{N}{3}\le a\le \frac{N}{2}$. Then $- \frac{N}{4}+\frac{3a}{4} \le \frac{N}{4}$; since $-\frac{N}{4}-\frac{a}{4}<-\frac{a}{2}$,  the assumption
(b) implies that
\begin{eqnarray*}
g(\frac{N}{4}-\frac{3a}{4}) - g(- \frac{N}{4}-\frac{a}{4}) \ge
g(\frac{N}{4}-\frac{3a}{4}) - g(-\frac{a}{2}) \geq
\int^{\frac{N}{4}-\frac{3a}{4}}_{-\frac{a}{2}} g^\prime(t)dt =
(*).\ens Since $-\frac{N}{4}\leq -\frac{a}{2}<
\frac{N}{4}-\frac{3a}{4}\leq 0$, (d)  implies that \bes (*)\geq
\int^{\frac{N}{4}-\frac{3a}{4}}_{-\frac{a}{2}}
g^\prime(-t-\frac{N}{2})dt =
\int^{-\frac{N}{2}+\frac{a}{2}}_{-\frac{3N}{4}+\frac{3a}{4}}
g^\prime(t)dt=(**).\ens
Since $-\frac{N}{2}+\frac{a}{2}\leq-\frac{N}{4}$ and
$g^\prime$ is increasing on $]-\infty,-\frac{N}{4}]$,
 shifting the
integration interval by $-\frac{N}{4}+\frac{3a}{4}\ge 0$ to the
left yields  that \bes (**)\geq
\int^{-\frac{N}{4}-\frac{a}{4}}_{-\frac{N}{2}} g^\prime(t)dt =
g(-\frac{N}{4}-\frac{a}{4});\ens thus \eqref{60208b} holds, as
desired. \ep

Proposition \ref{bn-4-1} immediately leads to the following simple criterion for
a function to belong to $\cap_{0<a<N}V_{N,a}.$
\bc \label{bn-19}
Let $N>0$, and
assume that a continuous and piecewise continuously differentiable function $g: \mr \to \mr$  with $\supp\ g =[-\frac{N}{2},\frac{N}{2}]$
 satisfies the following conditions:
\begin{itemize}
    \item[{\rm (a)}] $g$ is symmetric around the origin;
    \item[{\rm (b)}] $g^\prime$ is positive and increasing on $]-\frac{N}{2},0[$.
\end{itemize}
Then $g \in \cap_{0<a<N}V_{N,a}$.
\ec

We will now describe several functions belonging to $V_{N,a}$, either for
all $a\in ]0, N[$ or a subinterval hereof.

\bex \label{50911b} Consider the B-splines $B_N, \, N\in \mn,$ defined
recursively by \bes B_1 = \chi_{[-1/2, 1/2]}, \, B_{N+1} = B_N
*B_1.\ens
In  \cite[Lemma 1]{TS} it is proved that for  $N\in \mn\setminus\{1\}$,
\begin{itemize}
    \item[{\rm (i)}] $B_N^\prime$ is  increasing on $]-\frac{N}{2},-\frac{N}{4}+\frac{1}{4}]$;
    \item[{\rm (ii)}] $B_N^\prime (-x-\frac{N}{2}) \leq B_N^\prime (x)$ for $x\in [-\frac{N}{4},0[$.
\end{itemize}
Thus Proposition \ref{bn-4-1} implies that $B_N \in \bigcap_{0<a<N} V_{N,a}$
for $N\in \mn\setminus\{1\}$. \ep \enx

\bex \label{50911c} Let $N\in \mn\setminus\{1 \}$ and define
$$f_N(x):=\cos^{2N-2}(\frac{\pi x}{N})\chi_{[-\frac{N}{2},\frac{N}{2}]}(x).$$
Direct calculations show
that for $x\in ]-\frac{N}{2},-\frac{N}{4}]$,
$$f_N^{\prime \prime}(x)
=\frac{(2N-2)\pi^2}{N^2}
\cos^{2N-4}(\frac{\pi x}{N})
\left( (2N-3)\sin^2(\frac{\pi x}{N})-\cos^2(\frac{\pi x}{N})\right)\geq 0$$
and for $x\in [-\frac{N}{4},0[$,
$$f_N^\prime(x) - f_N^\prime(-x-\frac{N}{2})
=\frac{ (2N-2)\pi}{2N}\sin(\frac{2\pi x}{N})
\left(\sin^{2N-4}(\frac{\pi x}{N})-\cos^{2N-4}(\frac{\pi x}{N}) \right)\geq 0.$$
By Proposition \ref{bn-4-1}, $f_N \in\bigcap_{0<a<N} V_{N,a}$.
\ep \enx

In the following examples we consider continuous and compactly supported ``variants" of the two-sided exponential function, the Gaussian, and other classical functions.

\begin{ex} Let $N>0$ and define
$$h_N(x):=\left(e^{-|x|}-e^{-\frac{N}{2}}\right) \chi_{[-\frac{N}{2},\frac{N}{2}]}(x).$$
Then $h_N \in\bigcap_{0<a<N} V_{N,a}$ by Corollary \ref{bn-19}.
\ep
\end{ex}

\begin{ex} Let $N>0$ and define
$$k_N(x):=\left(\frac{1}{1+|x|}-\frac{1}{1+\frac{N}{2}}\right) \chi_{[-\frac{N}{2},\frac{N}{2}]}(x).$$
Then $k_N \in\bigcap_{0<a<N} V_{N,a}$  by Corollary \ref{bn-19}.
\ep
\end{ex}

In the following examples the simple sufficient conditions
in Proposition \ref{bn-4-1} and Corollary \ref{bn-19} are not satisfied. We will
use the definition directly to show that the considered functions belong
to $V_{N,a}$ for certain ranges of the parameter $a$.

\bex \label{50911e-1} Let $N>0$, and consider
\begin{equation*} 
p_N(x):=\left(\frac{1}{1+x^2}-\frac{1}{1+(\frac{N}{2})^2}\right)
\chi_{[-\frac{N}{2},\frac{N}{2}]}(x).
\end{equation*}
We will show that
\begin{itemize}
    \item[(a)] $p_N \in\bigcap_{\frac{3N}{7}\leq a<N} V_{N,a}$;
    \item[(b)] $p_N \in\bigcap_{\frac{N}{3}\leq a<\frac{3N}{7}} V_{N,a}$
    if $N\geq \sqrt{\frac{12}{5}}\approx 1.5451\cdots$.
\end{itemize} It is clear that (A1) and (A2) hold, so for the considered values for $a$
we now check (A3).  In fact the argument below will
prove more, namely that
\begin{equation}\label{bn-47}
\Delta_a^2 p_N(x) \geq 0, \ x\in
[-\frac{N}{2},-\frac{N}{4}+\frac{3a}{4}].
\end{equation} Considering any $a\in
[\frac{N}{3}, N[,$ we have
 $$-\frac{N}{2}<-\frac{N}{4}+\frac{3a}{4}<\frac{N}{2},\,
 a-\frac{N}{2}<-\frac{N}{4}+\frac{3a}{4}<a+\frac{N}{2}, \,
-\frac{N}{4}+\frac{3a}{4}<2a-\frac{N}{2};$$
so the fact that $p_N>0$ on
$]-\frac{N}{2},\frac{N}{2}[$, $p_N(\cdot-a)>0$ on
$]a-\frac{N}{2},a+\frac{N}{2}[$ and $p_N(\cdot-2a)>0$ on $]2a-\frac{N}{2},2a+\frac{N}{2}[$
immediately shows that
\begin{itemize}
   \item[(1)]  $p_N(x)\neq 0,\, x\in
]-\frac{N}{2}, -\frac{N}{4}+\frac{3a}{4}];$
    \item[(2)] $p_N(x-a)=0,
\, x\in [-\frac{N}{2},a-\frac{N}{2}]$;\,
$p_N(x-a)\neq 0, \, x\in ]a-\frac{N}{2}, -\frac{N}{4}+\frac{3a}{4}];$
    \item[(3)] $p_N(x-2a)=0, \, x\in
[-\frac{N}{2}, -\frac{N}{4}+\frac{3a}{4}].$
\end{itemize}
Thus
$$\Delta_a^2 p_N(x)= \left\{
\begin{array}{ll}
p_N(x), & x\in [-\frac{N}{2},-\frac{N}{2}+a]; \\
p_N(x)-2p_N(x-a), & x\in ]-\frac{N}{2}+a,-\frac{N}{4}+\frac{3a}{4}].
\end{array}
\right.
$$
Since $p_N(x)\geq 0, \ x \in [-\frac{N}{2},-\frac{N}{2}+a]$,
it is now enough to prove that
\begin{equation} \label{60212a}
p_N(x)-2p_N(x-a) \geq 0,\ x\in
]-\frac{N}{2}+a,-\frac{N}{4}+\frac{3a}{4}].
\end{equation}
Let $x\in ]-\frac{N}{2}+a,-\frac{N}{4}+\frac{3a}{4}]$. We see that
\begin{eqnarray} \notag
p_N(x)-2p_N(x-a) &=& p_N(x)-p_N(x-a)-p_N(x-a) \nonumber \\ \notag
&=& \frac{1}{1+x^2} -\frac{1}{1+(x-a)^2} -\left(\frac{1}{1+(x-a)^2}-\frac{1}{1+(\frac{N}{2})^2} \right) \nonumber \\ \notag
&=& \frac{(x-a)^2-x^2}{(1+x^2)(1+(x-a)^2)} -\frac{(\frac{N}{2})^2-(x-a)^2}{(1+(x-a)^2)(1+(\frac{N}{2})^2)}. \\ \label{bn-36}
\end{eqnarray}

In order to prove (a), we now assume that
 $\frac{3N}{7}\leq a<N$. Since
$x^2 \leq (\frac{N}{2})^2$
and
$(x-a)^2 -x^2 \geq 0$, we have
\begin{eqnarray*}
p_N(x)-2p_N(x-a)
& \geq &
\frac{(x-a)^2-x^2}{(1+(\frac{N}{2})^2)(1+(x-a)^2)} -\frac{(\frac{N}{2})^2-(x-a)^2}{(1+(x-a)^2)(1+(\frac{N}{2})^2)}\\
&=&
\frac{h(x)}{(1+(\frac{N}{2})^2)(1+(x-a)^2)},
\end{eqnarray*}
where
$$h(x):=2(x-a)^2-x^2-\frac{N^2}{4}=(x-2a)^2 -2a^2-\frac{N^2}{4}.$$
Note that the quadratic function $h$ is symmetric around
$x=2a.$ Since $ -\frac{N}{4}+\frac{3a}{4} <2a$ and $\frac{3N}{7}\leq a<N$,
we have
\begin{equation}\label{bn-35}
h(x)\geq h(-\frac{N}{4}+\frac{3a}{4})=\frac{1}{16}(N-a)(7a-3N) \geq 0.\
\end{equation}
Thus $$ p_N(x)-2p_N(x-a) \geq 0.$$ Therefore \eqref{bn-47} holds,
i.e., the proof of (a) is completed.

In order to prove (b), assume now that $N\geq \sqrt{\frac{12}{5}}$ and $\frac{N}{3}\leq a<\frac{3N}{7}$.
Since $(x-a)^2-x^2=2a(\frac{a}{2}-x)$ and
$(\frac{N}{2})^2-(x-a)^2=(\frac{N}{2}-x+a)(\frac{N}{2}+x-a)$, \eqref{bn-36} implies that
\begin{eqnarray*}
 p_N(x)-2p_N(x-a)
= \frac{1}{1+(x-a)^2}\left(
\frac{2a(\frac{a}{2}-x)}{1+x^2} -
\frac{ (\frac{N}{2}-x+a)(\frac{N}{2}+x-a)}{1+(\frac{N}{2})^2}\right).
\end{eqnarray*}
Using $-\frac{N}{2}+a<x\leq -\frac{N}{4}+\frac{3a}{4}$, it follows that
$\frac{N}{2}-x+a<N$ and $\frac{N}{2}+x-a\leq \frac{N}{4}-\frac{a}{4}\leq -x+\frac{a}{2}$; thus
\begin{eqnarray*}
 p_N(x)-2p_N(x-a)&\geq&
 \frac{1}{1+(x-a)^2}\left(
\frac{2a(\frac{a}{2}-x)}{1+x^2} -
\frac{ N(\frac{a}{2}-x)}{1+(\frac{N}{2})^2}\right) \\
&=&  \frac{\frac{a}{2}-x}{1+(x-a)^2}
\left(\frac{q(x)}{(1+x^2)( 1+(\frac{N}{2})^2 ) }\right),
\end{eqnarray*}
where $q(x):=2a(1+(\frac{N}{2})^2)-N(1+x^2)$.
Since $\frac{N}{3} \leq a <\frac{3N}{7}$ and $-\frac{N}{2}+a<x\leq -\frac{N}{4}+\frac{3a}{4}$,
we have
$-\frac{N}{6} < x\leq \frac{N}{14}$; so $x^2 < \frac{N^2}{36}$.
Thus
\begin{eqnarray*}
  q(x) \geq  2\frac{N}{3} \left(1+(\frac{N}{2})^2\right)-N\left(1+\frac{N^2}{36} \right) = N\left(-\frac{1}{3}+\frac{5N^2}{36}\right)
 \geq 0.
  \end{eqnarray*}
Thus $p_N(x)-2p_N(x-a) \geq 0.$
Therefore \eqref{bn-47} holds, as desired.\ep \enx


\bex \label{50911e} Let $N>0$, and consider
\begin{equation*} 
g_N(x):=\left(e^{-x^2}-e^{-\frac{N^2}{4}}\right)
\chi_{[-\frac{N}{2},\frac{N}{2}]}(x).
\end{equation*}
We will show that $g_N \in\bigcap_{\frac{3N}{7}\leq a<N} V_{N,a}$.
As in Example \ref{50911e-1}, see \eqref{60212a},
it suffices to check that
\begin{equation*} 
g_N(x)-2g_N(x-a) \geq 0,\ x\in ]-\frac{N}{2}+a,  -\frac{N}{4}+\frac{3a}{4} ]
\end{equation*}
Let $x\in  ]-\frac{N}{2}+a,-\frac{N}{4}+\frac{3a}{4}]$, and
$a\in
[\frac{3N}{7},N[$. We see that
\begin{eqnarray*}
g_N(x)-2g_N(x-a) &=& g_N(x)-g_N(x-a)-g_N(x-a) \\
&=&  e^{-x^2} - e^{-(x-a)^2} - \left( e^{-(x-a)^2}-e^{-\frac{N^2}{4}} \right) \\
&=&  e^{-(x-a)^2} \left( e^{(x-a)^2 -x^2}-1 \right)
- e^{-\frac{N^2}{4}} \left( e^{\frac{N^2}{4}-(x-a)^2}-1 \right).
\end{eqnarray*}
Since
$-(x-a)^2 \geq -\frac{N^2}{4}$
and
$(x-a)^2 -x^2 \geq 0$, we have
\begin{eqnarray*}
g_N(x)-2g_N(x-a)
& \geq &   e^{-\frac{N^2}{4}} \left( e^{(x-a)^2 -x^2}-1 \right)
- e^{-\frac{N^2}{4}} \left( e^{\frac{N^2}{4}-(x-a)^2}-1 \right)\\
&=&
e^{-\frac{N^2}{4}} \left( e^{(x-a)^2 -x^2}- e^{\frac{N^2}{4}-(x-a)^2}\right) \\
&=& e^{ -(x-a)^2}
\left( e^{h(x)}-1 \right),
\end{eqnarray*}
where
$h(x):=2(x-a)^2-x^2-(\frac{N}{2})^2.$
From \eqref{bn-35}, we have
$h(x)\geq 0$.
Thus $$ g_N(x)-2g_N(x-a) \geq 0,$$
as desired.\ep
 \enx

\noindent{\bf Competing interests:}
The authors declare that they
have no competing interests.

\noindent{\bf Authors' contributions:}
All authors contributed equally to this work. All authors read and approved the final manuscript.

\noindent{\bf Acknowledgments:} The authors would like to thank the reviewers for
many useful suggestions, which clearly improved the presentation. In particular, one
reviewer suggested to formulate the condition (A3) for membership of $V_{N,a}$
in terms of second order differences, which is much more transparent
than our original condition. He also suggested the approach in the current  proof of
Proposition \ref{bn-4-1}, which is shorter than our original proof.
This research was supported by Basic Science Research Program
through the National Research Foundation of Korea(NRF) funded by
the Ministry of Education (2013R1A1A2A10011922).

\end{document}